\documentclass[12pt]{article}
 
\usepackage[a4paper,margin=2cm]{geometry}

\usepackage{amssymb}
\usepackage{txfonts}
\usepackage[backend=biber,style=authoryear]{biblatex}
\usepackage[capitalise]{cleveref}
\usepackage{multirow}
\usepackage{float}
\usepackage{makecell}
\crefname{appendix}{}{}
\usepackage[figuresright]{rotating}
\usepackage{authblk}

\let\cite\parencite

\title{Recursive Target Body Approach for Low-Thrust Multiple Gravity-Assist Sequence Optimization}

\author[1]{Sean Cowan}
\author[1]{Ron Noomen}
\affil[1]{Delft University of Technology, Faculty of Aerospace Engineering, Kluyverweg 1, 2629HS, Delft, Netherlands}
\date{\today}

\begin{document}

\maketitle

\begin{abstract}

  This work aims to automate the design of Multiple Gravity-Assist (MGA) transfers between planets
  using low-thrust propulsion. In particular, during the preliminary design phase of space missions,
  the combinatorial complexity of MGA sequencing is very large, and current optimization approaches
  require extensive experience and can take many days to simulate. Therefore, a novel optimization
  approach is developed here -- called the Recursive Target Body Approach (RTBA) -- that uses the
  hodographic-shaping low-thrust trajectory representation together with a unique combination of
  tree-search methods to automate the optimization of MGA sequences. The approach gradually
  constructs the optimal MGA sequence by recursively evaluating the optimality of subsequent
  gravity-assist targets. Another significant contribution to the novelty of this work is the use of
  parallelization in an original way involving the Generalized Island Model (GIM) that enables the use
  of new figures of merit to further increase the robustness and accelerate the convergence. An
  Earth-Jupiter transfer with a maximum of three gravity assists is considered as a reference
  problem. The RTBA takes 21.5 hours to find an EMJ transfer with 15.4 km/s $\Delta V$ to be the
  optimum. Extensive tuning improved the quality of the MGA trajectories substantially, and as a
  result a robust low-thrust trajectory optimization could be ensured. A distinct group of highly
  fit MGA sequences is consistently found that can be passed on to a higher-fidelity method. In
  conclusion, the RTBA can automatically and reliably be used for the preliminary optimization of
  low-thrust MGA trajectories. 

\end{abstract}



\section{Introduction}




Large-scale interplanetary space missions, as of 2025, are extending humanity’s reach to the
periphery of the Solar System. These missions are instrumental in advancing our understanding of the
universe and assessing the potential for life beyond Earth (e.g., JUICE \cite{juicemission}, LISA
\cite{lisamission}, Europa Clipper \cite{clippermission}). However, the challenges associated with
such missions are substantial, encompassing long duration, high cost, and significant
technological and engineering complexity. Among these challenges, trajectory design plays a
critical role.

Current chemical propulsion systems, while well-established, are insufficient for directly reaching
the outer regions of the Solar System without external assistance. To overcome this limitation,
Gravity-Assist (GA) maneuvers have been developed, enabling spacecraft to gain energy through
planetary flybys. Simultaneously, efficient low-thrust propulsion technologies, such as ion and
Hall-effect thrusters, have emerged as viable alternatives, offering high exhaust velocities and
much improved propellant efficiency. Notable missions like BepiColombo and DART have demonstrated the
potential of low-thrust engines, providing a foundation for future exploration. Despite these
advancements, the integration of low-thrust technologies and GA maneuvers into mission design presents
formidable optimization challenges, particularly in the context of designing Multiple Gravity-Assist
(MGA) sequences.

Many methods exist to solve and optimize low-thrust optimization problems. Examples include but are
not limited to those involving Pontryagin's minimum principle \cite{pontryagin2018mathematical,
rosa2006genetic, zhang2015low} and those involving the direct optimization of a Non-Linear
Programming Problem (NLP) \cite{vasile2002design, carnelli2009evolutionary, morante2020hybrid}. In
the direct optimization methods, shooting methods are commonly used. However, since these are
computationally expensive, numerous (semi-)analytical alternatives have been developed over the past
decades. A widely spread approach is the traditional combination of Lambert arcs together with
impulsive $\Delta V$'s that can be split into smaller and smaller segments, thereby more closely
approximating a continuous thrust trajectory -- coined as Sims-Flanagan transcription
\cite{sims1999preliminary}. Another category of analytical approaches uses averaging
\cite{gao2007near, kechichian2009optimal}. A completely different approach is to use semi-analytical
shaping methods. Examples of shaping methods are exponential sinusoids \cite{petropoulos2004shape},
inverse polynomials \cite{wall2009shape}, spherical shaping \cite{novak2011improved}, and Fourier
series \cite{abdelkhalik2012dynamic}. 

The aspect of MGA sequence design has historically relied on heuristic methods, such as Tisserand
graphs and experience-driven initial guesses \cite{crain2000interplanetary, debban2002design,
strange2002graphical}. However, modern mission design increasingly leverages computational
optimization and automation, which provide greater flexibility and efficiency. Optimization methods
used can be split into single-loop and nested-loop approaches. For the former, the design variable
of the various flybys has to be encoded in a single vector together with the variables relating to
the transfer leg \cite{gad2011hidden, englander2012automated, chilan2013automated,
abdelkhalik2012dynamic}. The latter approach involves using various combinations of heuristic
optimization techniques like genetic algorithms or ant colony optimization as well as tree-search
methods such as beam-search \cite{englander2012automated, englander2017automated,
chilan2013automated, CERIOTTI20101202}. Despite the progress, current optimization algorithms face
significant limitations, including high computational runtimes \cite{englander2017automated} and
limited validation capabilities. Moreover, the coupling of MGA sequencing with low-thrust trajectory
optimization (LTTO) exacerbates the complexity of the problem due to the combinatorial nature of
sequence design and the multi-dimensional nature of low-thrust trajectory modeling. Comparisons of
single- and nested-loop approaches lean toward the use of nested-loop approaches
\cite{ellison2018robust, morante2021survey}, which is the approach that is used in this study. When
applied to asteroid targeting though, a mixed-code single-loop approach was shown to be better
\cite{zhang2015analysis}. It clearly is an area with potential for future research.

Many of the aforementioned studies have focused on impulsive $\Delta V$ trajectories, and only
recently a number of papers have emerged that also incorporate low-thrust propulsion into the
dynamics description \cite{englander2017automated, morante2019multi, fan2021fast, fan2022improved}.

Separately, a subset of studies that tackle the low-thrust MGA optimization problem uses tree-search
methods \cite{ellison2018robust, CERIOTTI20101202, bellomemulti, vasile2015incremental,
hennes2015interplanetary, james2017ananalysis}. These methods have the advantage that they are
discrete by nature and therefore serve well for combinatorial problems.

This study addresses this computational challenge by proposing a novel approach to optimize
low-thrust MGA trajectories. The method employs a nested-loop framework, wherein the inner loop
performs LTTO for specific MGA sequences, and the outer loop explores the combinatorial search space
to identify optimal sequences. A key innovation is the combination of tree-search methods to enhance
convergence. To further improve computational
efficiency and confidence in sequence evaluation, the Generalized Island Model (GIM) is adapted to
enable parallel evaluations. This approach integrates hodographic-shaping methods for trajectory
approximation, balancing computational cost and fidelity in ranking MGA sequences. Together, these
techniques are coined the Recursive Target Body Approach (RTBA).

Additionally, this work introduces a new metric for quantifying sequence optimality, providing a
more informed basis for decision-making in mission design. By addressing the interplay between
sequence optimization and LTTO fidelity, the proposed framework aims to reduce computational
overhead while maintaining or improving solution quality. These advancements hold the potential to
significantly enhance the performance and reliability of trajectory optimization processes for
low-thrust interplanetary missions.

The paper is organized as follows. \cref{sec:theorydescrip} lays out the theoretical description of
the dynamics as well as the techniques employed in the RTBA. \cref{sec:optimapproach} introduces the
optimization approach, presenting both the structure as well as the tuning. \cref{sec:results}
subsequently presents the main results. \cref{sec:conclusion} draws conclusions on the work done in
this study.

\section{Theoretical description} \label{sec:theorydescrip}

\subsection{Dynamics} \label{sec:dynamical-description}

To describe a dynamics system, a reference system is required. For the heliocentric trajectories,
the ECLIPJ2000 reference frame \cite{eclipj2000} is used due to its ecliptic nature and widespread use. For the
description of gravity-assist maneuvers, two frames are needed. First, a TNW-frame is taken as an
inertial reference frame centered at the Solar System Barycenter (SSB). This frame is useful when
using a patched-conics approximation, as it simplifies the expressions. Second, a separate frame is
needed to describe the outgoing velocity vector relative to the incoming one, which will be denoted
as the 'local' frame. To accompany the reference frames, a cylindrical coordinate system is chosen,
which is useful for the trajectory description \cite{gondelach2015hodographic}. Since this work aims
to develop a preliminary optimization method for interplanetary missions as a primary goal, two-body
dynamics are used with a continuous perturbative force for the thrust. A more extensive description
is given in \cref{sec:appendix}.

\subsection{Hodographic shaping} \label{sec:hodo-shaping}

The dynamics description that is used to establish multiple gravity-assist trajectories is now
complete. However, from \cref{equation:twobodymotionpert}, a solution still needs to be found. For
this, a shaping method is used, due to its computational efficiency and relatively high
accuracy \cite{petropoulos2004shape, wall2009shape, novak2011improved, abdelkhalik2012dynamic}.
Specifically, the hodographic-shaping method is used \cite{gondelach2015hodographic}. This method
combines a number of elements that most shaping methods lack: a direct time-of-flight solution, and
a description that allows for three dimensions, many revolutions, and a thrust acceleration limit
through constraint optimization.

Hodographic shaping relies on the velocity hodograph, which plots velocity components as functions
of time or polar angle; the time-based approach is preferred for its physically real time-derivative
and simpler implementation. The functions are in cylindrical coordinates as mentioned in
\cref{sec:dynamical-description}.

Velocity-hodograph shapes require mathematical descriptions. These functions are built from
analytically differentiable and integrable base functions, such as a sine, cosine, exponential,
power series, or combinations thereof. A minimum of three base functions per axis ensures boundary
condition satisfaction, resulting in three coefficients per axis. Additional velocity functions and
coefficients can be inserted to enhance the method's adaptability to describe lower $\Delta V$
solutions. The recommended base and additional functions are given below
\cite{gondelach2015hodographic} and are based on an Earth-Mars transfer. The optimal set of
functions for other transfers may be different, but the characteristics of a transfer to any outer
planet are mostly the same and for a transfer to an inner planet the drop in performance of these
functions was small when compared to the inner-planet recommended functions. \cref{eq:vbasern} and
\cref{eq:vbasea} are the base functions that are required for any hodographically-shaped trajectory.
Three coefficients can be seen per axis, which allow the boundary conditions to be met.
\cref{eq:vaddrn} and \cref{eq:vadda} are the additional velocity functions that one can add to the
radial, normal, and/or axial component of the velocity; each extra degree of freedom adds three
free coefficients (one per axis).

\begin{equation} \label{eq:vbasern}
      V^{base}_{r,n} = c_1+c_2t+c_3 t^2 \\
\end{equation}
\begin{equation} \label{eq:vbasea}
      V^{base}_{a} = c_4 \cos(2\pi t (N+0.5))  +c_5 t^3 \cos (2\pi t (N+0.5))  +c_6 t^3 \sin(2 \pi t (N + 0.5)) \\
\end{equation}
\begin{equation} \label{eq:vaddrn}
      V^{add}_{r,n} = c_7 t  \sin(0.5 t \pi)  +c_8t  \cos(0.5 t\pi) + c_9 ...\\
\end{equation}
\begin{equation} \label{eq:vadda}
  V^{add}_{a} = c_{10} t^4 \cos(2 t \pi (N + 0.5)) +c_{11} t^4 \sin(2 t \pi (N + 0.5)) + c_{12} ...
\end{equation}

From the velocity functions, the analytical integrals are calculated to have an expression for the
change in position. The derivatives are calculated to determine the inertial acceleration. Using
\cref{equation:twobodymotionpert}, the thrust acceleration can be derived. The thrust acceleration
is then integrated over time to produce the total $\Delta V$ required for that particular
trajectory.

\section{Optimization approach} \label{sec:optimapproach}

\subsection{The inner loop} \label{subsec:ltto}

The low-thrust trajectory optimization (LTTO) problem is considered as a Hybrid Optimal Control
Problem (HOCP) and formulated as a Mixed-Integer Non-Linear Program (MINLP) using hodographic
shaping, discussed in \cref{sec:hodo-shaping}. Hodographic shaping contains both continuous variables
(e.g. coefficients of a shaping function), as well as discrete ones (e.g. number of revolutions).


For the optimization itself, PyGMO \cite{biscani2020parallel} was used. The Simple Genetic Algorithm
(SGA) is chosen due to its relative simplicity in tuning and its widespread use. The objective
function for the LTTO will only consider $\Delta V$, but alternative objectives like delivery mass,
mass fraction, and propellant mass can easily be tested as well. 

The input parameters are the MGA sequence, the initial mass and specific impulse of the spacecraft,
the departure and arrival conditions of the first and final target body, the number of free
parameters -- the number of added base functions to the velocity functions -- and finally the 
generation count and population size of the genetic algorithm.

As part of the design space exploration, to understand and improve the results of each individual
low-thrust trajectory, numerous tuning steps were conducted. Besides tuning which design variables
were best included in the optimization, the size of their bounds and the genetic algorithm parameters
were also investigated.

The design variables -- shown in \cref{table:lttovariables} -- are split into three categories:
general, leg-specific, and GA-specific parameters. The bounds that were found to be optimal are also
given. A departure window of 60 days was used as it allowed for a more exhaustive search of that window
given a set number of islands that search. To ease the search space, some departure and arrival
conditions are set to 0, which means the spacecraft departs/arrives at a state equal to that of the
departure/arrival planet. The time-of-flight bound encompassed almost all values that were found for
an Earth-Earth-Earth-Mars-Jupiter (EEEMJ) transfer. Similarly, the rest of the leg and GA-specific
parameters were determined based on the range of values observed for an EEEMJ transfer.

\begin{table}
    \caption{Design variables for LTTO optimization. Subscripts $g$ and $i$ refer to parameters that
    have multiple variables equal to the number of GAs or legs, respectively. For the departure
    date, any value can be chosen.}
    \begin{center}
        \begin{tabular}{c|ccc|c}
            \hline
            Variable Type & Design Variable Name & Variable & Unit & Value/Bound \\ \hline
            \multirow{7}{*}{General}& Departure Date & $t_{dep}$ & days & $[t_{dep}, t_{dep} + 60]$ \\
                                    & Arrival Velocity & $v_{arr}$ & m/s & $0$\\
                                    & Departure Velocity & $v_{dep}$ & m/s & $0$ \\
                                    & Departure In-plane Angle & $\theta_{dep}$ & rad & 0 \\
                                    & Departure Out-of-plane Angle & $\phi_{dep}$ & rad & 0\\
                                    & Arrival In-plane Angle & $\theta_{arr}$ & rad & 0 \\
                                    & Arrival Out-of-plane  & $\phi_{arr}$ & rad & 0 \\
            \hline
            \multirow{3}{*}{Leg specific} & ToF & $tof_i$ & days & $[100, 4500]$ \\
                                          & Free Velocity Coefficients & $(c_{10})_i, ..., (c_j)_i$
                                          & - & $[-3 \cdot 10^4, 3 \cdot 10^4]$ \\
                                          & Revolutions & $rev_i$ & - & $[0, 2]$\\
            \hline
                                          &  Incoming Velocity & $v_{g}$ & m/s & $[0, 5000]$ \\
                                          & GA altitude above planetary surface & $h_{p, g}$ & m & $[ 2 \cdot 10^5, 5 \cdot 10^{10}]$ \\
                                          & Orbit orientation angle & $\beta_g$ & rad & $[0, 2\pi]$\\
                                          & GA In-plane angle & $\theta_{g}$ & rad & $[0, 2\pi]$\\
            \multirow{-6}{*}{GA specific}& GA Out-of-plane angle & $\phi_{g}$ & rad & $[-\frac{1}{2}\pi, \frac{1}{2}\pi]$\\
            \hline
        \end{tabular}
    \end{center}
    \label{table:lttovariables}
\end{table}

The parameters related to the optimization algorithm are shown in \cref{table:generaltuning} in
\cref{sec:results}. These values were found based on the convergence and expected optimality of the
aforementioned EEEMJ transfer. Note that the number of generations is fixed to a value of 300; no
explicit convergence test nor criterion needs to be included. The number 300 was found to guarantee
good, converged solutions for all cases studied here. Other values resulted in premature
convergence, or a larger search space in the case of the free parameter count which prevented the
hodographic shapes from being optimized. The free parameter count represents the number of extra
base functions that are added to each axis of the shaping function. The last parameter 'Topology
probability' is related to a key part of the inner loop which is why the GIM is used. The GIM allows
for an accelerated optimization, by allowing multiple optimization processes to be run in parallel
that can exchange high-fitness value individuals in a population between islands. In this case, it
is chosen to connect every island with every other island that has the same Target Body (TB). This
term is explained further in \cref{subsec:mgaso}. 

Without going into detail and to summarise how the main tuning process improved the quality of the
results, the difference in $\Delta V$ is tabulated from before and after the tuning process in
\cref{table:dvcompltto}, and then compared to \cite{fan2021fast}. For the MGA sequences EEMJ and
EEEMJ, the tuned configuration also returned competitive results. 

\begin{table}
  \caption{Initial comparison of minimum $\Delta V$ for the Earth-Jupiter and Earth-Mars-Jupiter
  transfer with results from \cite{fan2021fast}.}
  \label{table:dvcompltto}
  \begin{center}
    \begin{tabular}[c]{l|cc}
      \hline
      & \multicolumn{2}{c}{ Transfer} \\ \cline{2-3}
      \multirow{-2}{*}{} & EJ [km/s] &  EMJ [km/s] \\
      \hline
      Before tuning & 17.78 & 33.03 \\
      After tuning & 17.48 & 15.30 \\
      \cite{fan2021fast} & 17.81 & 15.18 \\
      \hline
    \end{tabular}
  \end{center}
\end{table}

\begin{figure}
    \centering
    \includegraphics[width=0.85\textwidth]{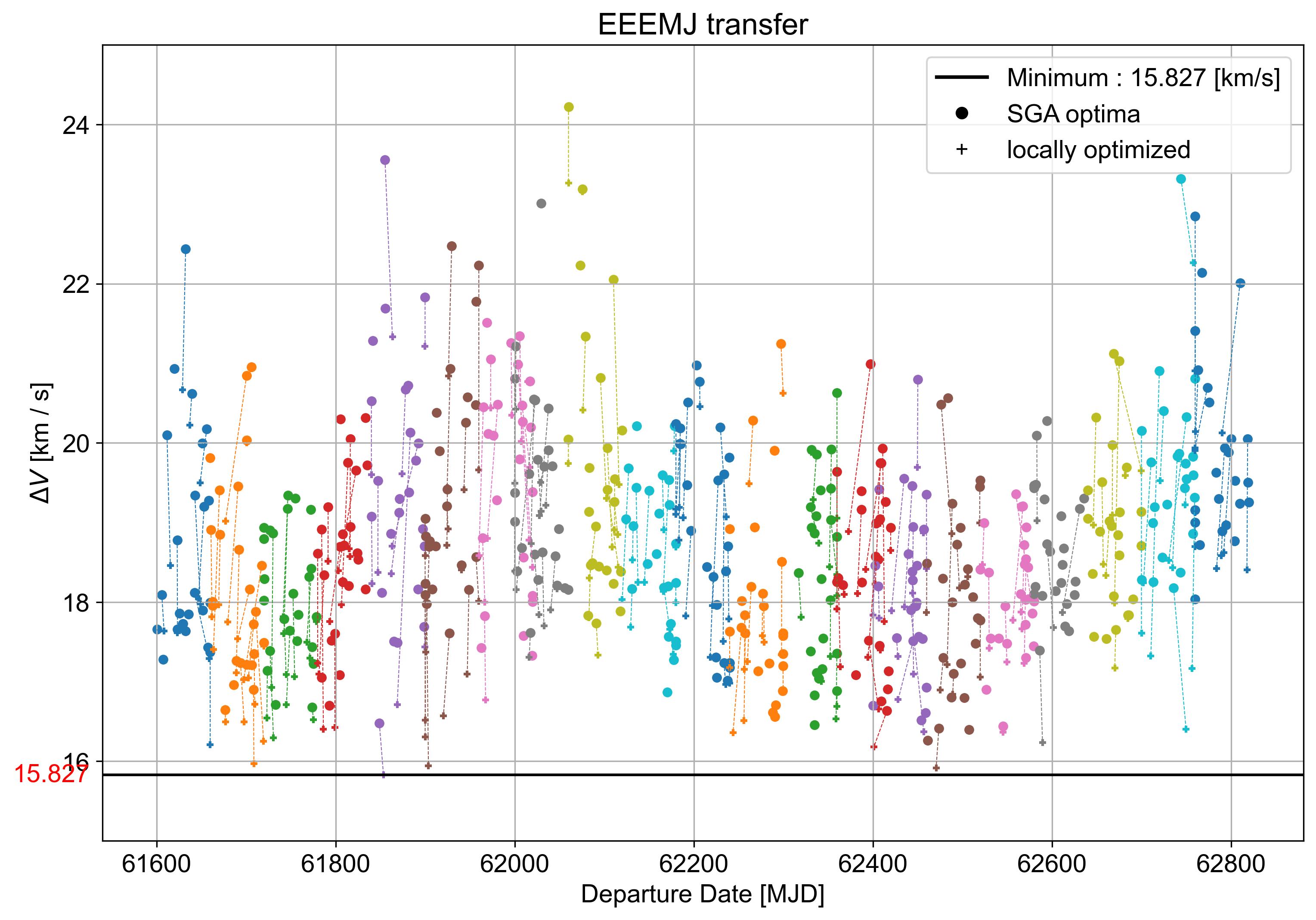}
    \caption{Grid search of departure dates from 61400-62860 MJD with and without additional local optimization step.}
    \label{figure:gsddate60localopt}
\end{figure}

To verify the performance and behavior of the optimization, a grid search was performed on an EEEMJ
transfer, the results for which are shown in \cref{figure:gsddate60localopt}, which can be compared
to Figure 12 of \cite{fan2021fast}. EEEMJ was chosen because it has three GAs and therefore
represents one of the more challenging transfers considered in this paper. In
\cref{figure:gsddate60localopt}, the grid steps are represented by 24 islands that were used to
optimize 60 day time intervals -- where each color represents a separate departure date interval and
each circle the optimum for a particular island. To verify the convergence of the optimization a
local gradient-based optimization step is added (represented by the '+' symbols). For the local
optimization, a Nelder-Mead Simplex algorithm is used. It can be seen that the transfers found with
the local optimization are almost always more optimal, sometimes with a significant difference in
$\Delta V$ and sometimes with barely any difference. This is caused by the random nature of the SGA,
which may not necessarily converge to the local minimum. The variation in $\Delta V$ across
departure dates can be explained by the optimality of the relative phasing of Earth, Mars, and
Jupiter; this variation is also confirmed by \cite{fan2021fast}. At a glance, it seems that the
local optimization causes different results, because the optimal departure date with only SGA
optimization is very different from the optimal departure date including local optimization.
However, if multiple SGA-optimized islands are compared around both 61800 MJD and 62400 MJD, as well
as locally optimized islands around the same departure date, then their $\Delta V$ values differ
only slightly. Between the most optimal SGA sequence and most optimal locally-optimized island,
there is a difference of two km/s. As will be seen in \cref{sec:results}, the pool of low-$\Delta V$
sequences differ by a margin of more than two km/s, and so therefore this difference is not
crucial. The quality and robustness of solutions is more than sufficient for the purpose of this new
RTBA tool.

\subsection{The outer loop} \label{subsec:mgaso}

As was established before, the outer loop is constructed using tree-search methods, as they are
considered most reliable when optimizing an inherently categorical search space
\cite{ellison2018robust, CERIOTTI20101202, bellomemulti, vasile2015incremental,
hennes2015interplanetary, james2017ananalysis}. Two main contenders for performant tree-search
methods that have been used to optimize the MGA sequencing problem consider almost exclusively
high-thrust cases: \cite{hennes2015interplanetary} and \cite{ellison2018robust}. The former
implemented the Beam Search Strategy (BSS) which is a greedy approach. The latter implemented
Monte-Carlo Tree Search (MCTS), which includes multiple randomized iterations and a gradual build-up
to the optimum.

\begin{figure}
  \begin{center}
    \includegraphics[width=0.75\textwidth]{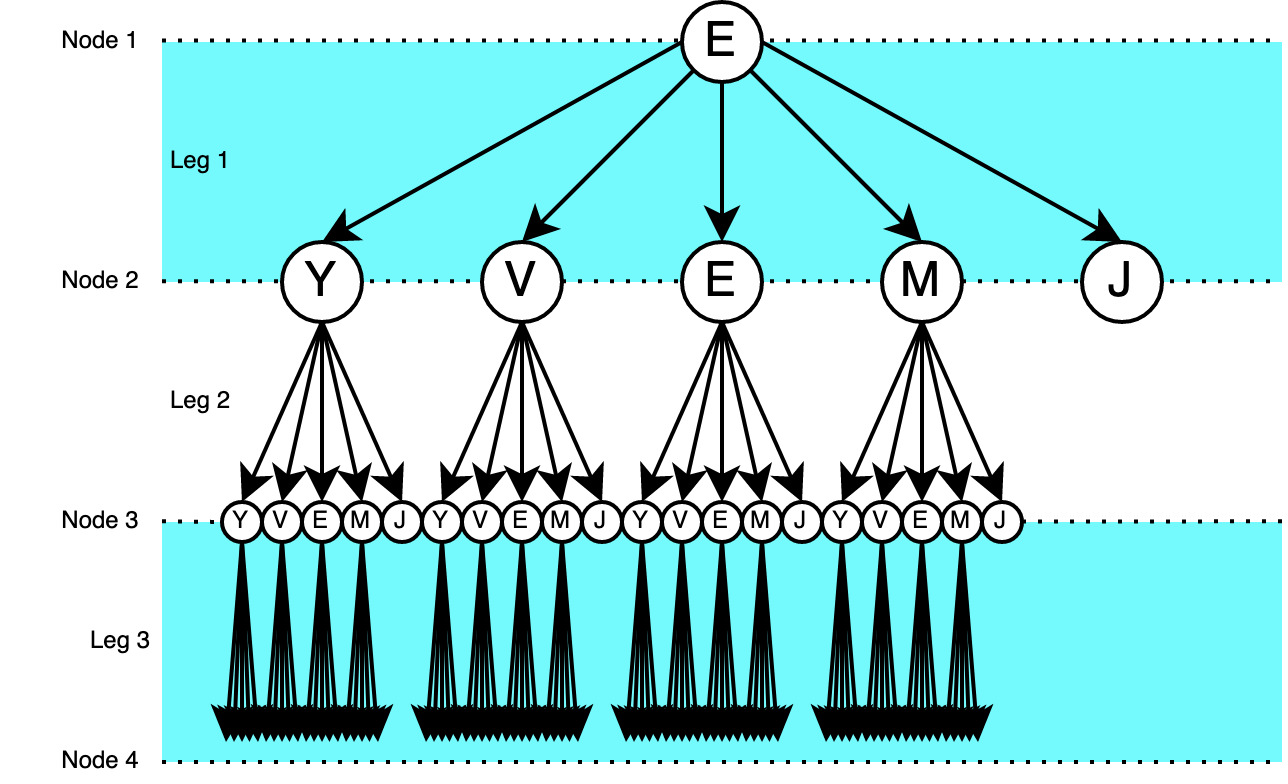}
  \end{center}
  \caption{A three-level tree for an Earth-Jupiter transfer. 'Y' is Mercury, and other planets are
    denoted by the first letter of their names.}
  \label{figure:combtree}
\end{figure}

In this work the BSS and MCTS are combined, which is a very novel element of this study, and this
combination is shown to produce a robust, greedy method that can prune possibilities effectively.
While MCTS is generally considered an effective tree-search method, given the run-time needed to
optimize a single sequence with the methodology described in this section, a greedy derivative is
chosen to speed up the optimization. The combinatorial space of all possible MGA sequences can be
represented as a tree that is directed, acyclic, and unweighted, as seen in \cref{figure:combtree}.
It is directed to visualize complete combinatorial enumeration, acyclic to avoid having to assess
identical $\Delta V$ costs for repeating sequences, and unweighted because the $\Delta V$ values
between two nodes are not always the same, but rather they depend on the specific and complete
sequence; $\Delta V$ values are taken as they are, without scaling and/or an extra penalty (for
constraint violation).

The RTBA is a greedy approach applied to low-thrust MGA sequencing. The representation of low-thrust
trajectories involves significantly more parameters than any high-thrust equivalent, resulting in
the complex and challenging optimization problem considered in this work. The RTBA consists of a few
steps in every recursion: Monte Carlo search, LTTO optimization, and Evaluation and pruning.

The RTBA is a recursive and greedy algorithm that repeats the same steps until specific criteria are
met. A few terms need to be introduced: First, the Pseudo Sequence (PS) is defined as the MGA
sequence that is gradually fixed after each subsequent recursion of the RTBA, and second a Target
Body (TB) is a GA candidate that immediately follows the fixed PS in each recursion.

Each recursion starts with a Monte-Carlo search, where random sequences are selected and evaluated.
The sequences are not entirely random, but rather evenly distributed across the TBs. Sequences are
then optimized with the LTTO process described in \cref{subsec:ltto}, and once the best TB has been
determined in that recursion, it is fixed in the PS and a new recursion is started, only considering
the sub-tree starting from that new PS. The RTBA uses recursion until there are no branches left, or
until the maximum allowed number of recursions has occurred. In the end, a ranking is given of all
sequences that were found, so the final optimal sequence is not guaranteed to be the PS, or have
a sequence length equal to the number of allowed GAs.

An extra element of novelty is the way in which the objective is considered at the sequence level,
and also at the TB level in each recursion. Rather than only using the minimum or mean $\Delta V$, a
linear combination is used. For the former case, it is used to assess the fitness of a single
sequence using the $\Delta V$'s of all islands ($f_s$). For the latter case, the fitness of TB
candidates -- meaning all sequences that use a given TB in a particular recursion -- are evaluated
using this linear combination of all sequences that flyby that TB ($f_{TB}$). In general, this helps
identify optimal trajectories by considering a wide search space with outliers, while also not
prematurely converging onto it. The ratio of the minimum and mean $\Delta V$ at the sequence level
and at the TB level is defined as the sequence fitness proportion ($\xi$) and the target body
fitness proportion ($\chi$). For example, if $\chi$ is 1.0, then only the sequence with the minimum
$\Delta V$ is used to evaluate the fitness of a given sequence ($f_s$). If $\chi$ is 0.7, then the
minimum $\Delta V$ contributes 70\% to the optimality of the sequence and the mean $\Delta V$
contributes 30\%. Analogously, if $\xi$ is 1.0, then only the sequence with the minimum $f_s$ is
used to evaluate the fitness of a given TB. This gives a useful, extra element of versatility to the
RTBA.

As the GIM is used, the number of CPUs needs to be specified, accompanied by a parameter indicating
the number of islands that need to be evaluated per sequence, denoted as $p$. Moreover, to quantify
the number of sequences that are evaluated, a number of parameters must be defined related to the
combinatorial complexity of the problem. The combinatorial complexity of any tree is defined as the
number of possible combinations, given by \cref{eq:frac}. This definition can be used at any
recursion level for every possible tree and sub-tree.

\begin{equation}
  C = \sum^n_{i=0} m^i
  \label{eq:frac}
\end{equation}

In \cref{eq:frac}, $C$ is the total number of combinations, $m$ is the number of possible GA
candidates, and $n$ is the maximum number of GAs remaining. If three GAs are allowed, and there are
eight possible GA planets, \cref{eq:frac} returns a complexity of $8^3 + 8^2 + 8^1 + 8^0 = 585$ unique sequences in the first
recursion. With the complexity of the problem defined, the proportion of the combinatorial space
that is evaluated during the optimization needs to be established. Parameter $q$ is therefore
defined to be the fraction of the combinatorial complexity $C$ that is evaluated during a single
recursion. For example, if $q$ is 0.3, then during that recursion the total remaining complexity $C$
is calculated and $0.3 \cdot C$ sequences are evaluated (this value is rounded up to the nearest
integer). For this work, $q$ is constant across the recursion levels, giving a total evaluated
fraction of $Q$.

A few tricks are implemented to augment the MGASO process. First and foremost, the possible GA
candidates are limited to the planets with a semi-major axis equal or smaller than the arrival
planet. This prunes Earth-Jupiter transfers like ENYJ (where Y is Mercury), which are likely not
optimal. Second, a feature present in most MCTS implementations is that every sequence is evaluated
only once; subsequent recursions cannot evaluate the same sequence. Third, sequences that flew by a
given TB will be tracked throughout all recursions and used in the calculation of the fitness of
that TB in any given recursion.

After tuning $p$, it was found that -- given the maximum disposable number of CPUs of 48 -- $p$
should be 14, using 42 CPUs and therefore evaluating exactly three sequences at each moment. This
allows for more precise estimation of the run time based on the combinatorial complexity.


\section{Results} \label{sec:results}


\begin{table}
  \caption{Tuned problem parameters for LTTO and MGASO loop.}
  \label{table:generaltuning}
  \begin{center}
    \begin{tabular}{c|c}
      \hline
      LTTO Parameter & Value \\ \hline
      Population size & 1200 \\
      Generation count & 300 \\
      Free parameter count & 1 \\
      Topology probability & 0.01 \\
      \hline\hline
      MGASO Parameter & Value \\ \hline
      $p$ & 14 \\
      \#CPUs & 42 \\ 
      Recursion count & 2 \\
      \hline
    \end{tabular}
  \end{center}
\end{table}

\begin{figure}
  \centering
  \includegraphics[width=\textwidth]{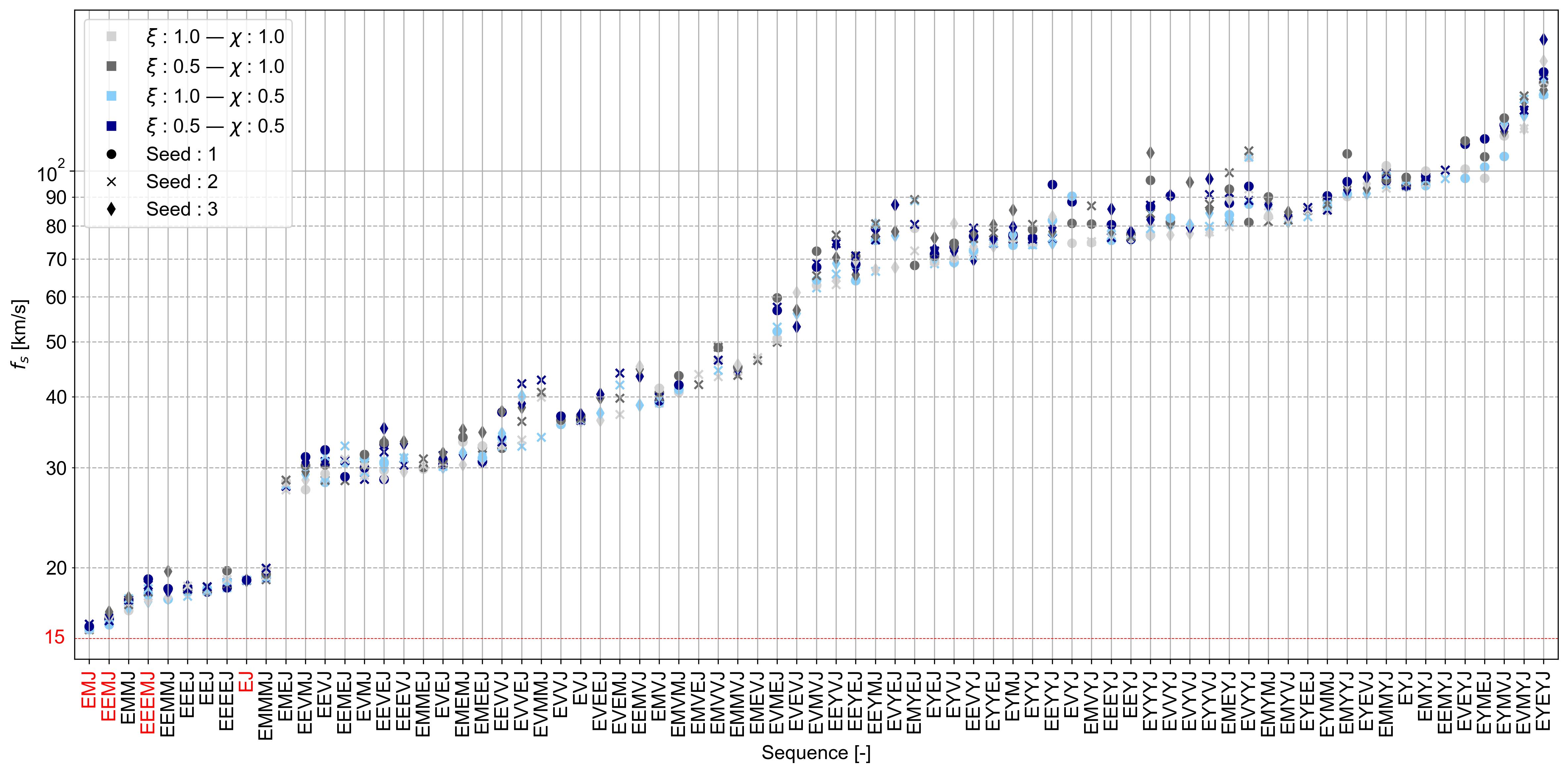}
  \caption{MGA sequences sorted by $f_s$, for $q = 0.5$ for three different seeds and four
  combinations of $\xi$ and $\chi$. Red sequences represent those present in \cite{fan2021fast}.}
  \label{figure:results}
\end{figure}

The main test case is an Earth-Jupiter transfer, which was also used in \cite{fan2021fast}. The
maximum number of GA's is set to three as it allows for a reasonably high combinatorial complexity
and it is comparable to literature \cite{englander2017automated, morante2019multi, fan2022improved}.
Specifically, the design variable bounds are given by \cref{table:lttovariables}, where the
departure date is specified to be between 61400 and 63400 MJD (or about 5.5 years). The MGASO
configuration is summarized in \cref{table:generaltuning}. For these simulations, 42 CPUs were used
on a single Intel Xeon Compute Node of the DelftBlue supercomputer, and $p$ (the number of islands
per sequence) is set to 14, as mentioned before. The maximum permissible number of recursions is
chosen to be two, because the test case allows for three GAs and therefore the final step will
already have fixed two of the three possible GA candidates. With these settings the MGASO was run
for $q=0.5$ (and consequently $Q=0.58$), the resulting sequences and their corresponding optimality
in terms of $f_s$ are shown in \cref{figure:results}. $q$ values of 0.1 and 0.3 were also tested but
are not shown for brevity, though the observed trends are identical, generally including a
significant portion of the most optimal sequences despite the smaller number of sequences evaluated.
Higher values of $q$ would start approaching an exhaustive search of all permissible transfers,
which would defeat the purpose of the RTBA and are therefore not tested.

\begin{figure}
  \centering
  \includegraphics[width=0.5\textwidth]{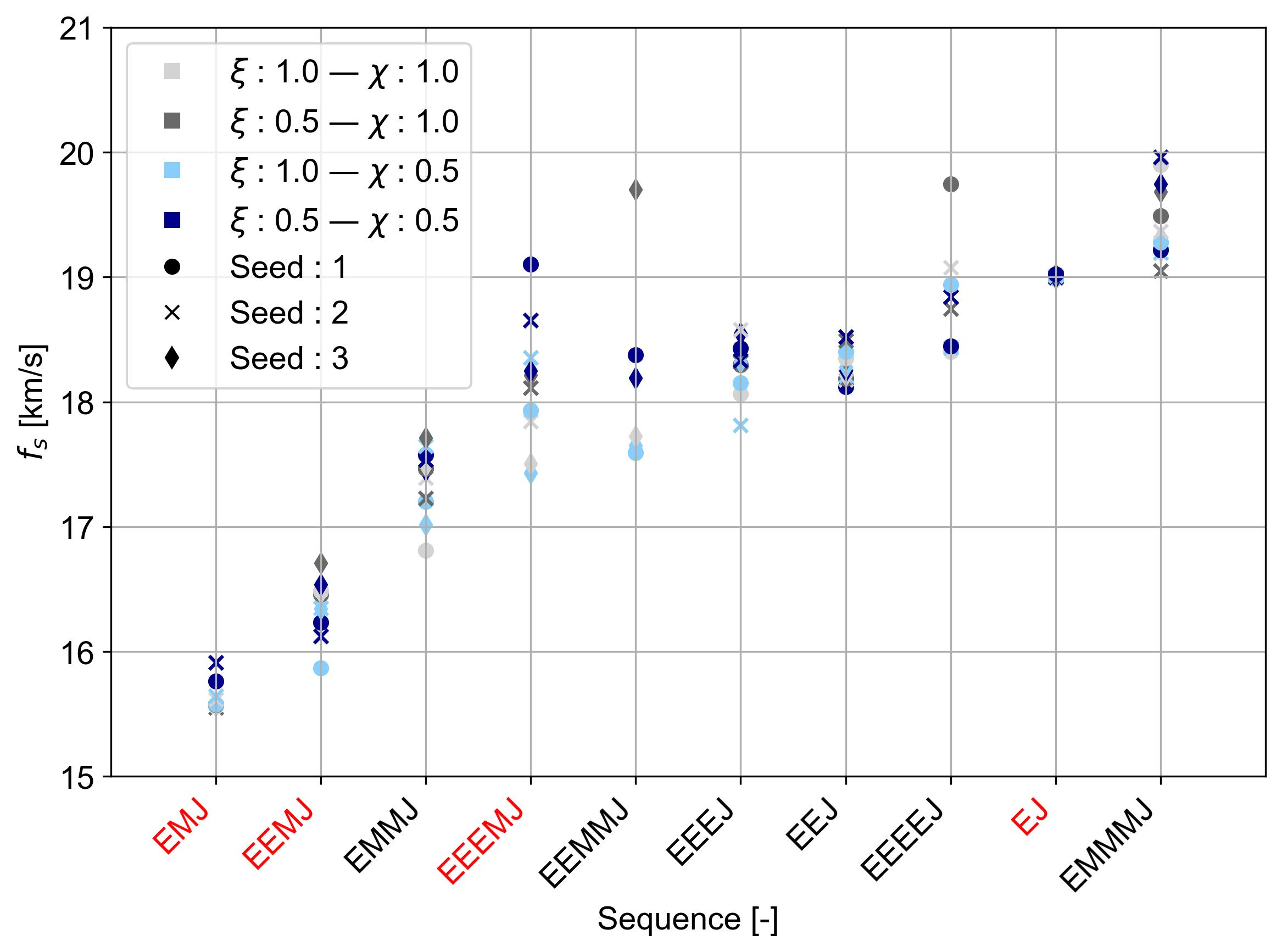}
  \caption{Most optimal MGA sequences sorted by $f_s$, for $q = 0.5$ for three different seeds and four
  combinations of $\xi$ and $\chi$. Red sequences represent those present in \cite{fan2021fast}.}
  \label{figure:results_highfit}
\end{figure}

The minimum $f_s$ observed is the EMJ transfer at about 15.4 km/s across all $q$ values, consistent with
\cref{table:dvcompltto}. This $f_s$ remains stable across $\xi$, $\chi$, and $q$, indicating a small
difference between the minimum and mean $\Delta V$, and similar findings across islands.

The spread of $f_s$ values (the difference between maximum and minimum $f_s$) increases with higher
$f_s$ values (mind the logarithmic scale). It ranges from about 1 km/s for low-$f_s$ sequences to 10-30
km/s for high-$f_s$ sequences. This trend stems from LTTO optimization tuning for specific
transfers; sequences diverging significantly from the transfers used for tuning show weaker convergence.

A major jump in $f_s$ (about 10 km/s) appears across all tested $q$ values. Sequences below this
threshold predominantly use sub-sequences ('MM', 'EM', 'EE') that are more likely to be optimal from
a theoretical point of view, while those above contain sub-optimal sequences (e.g., 'ME') or target
Venus or Mercury. The sub-optimality of these planets as GA targets is in line with astrodynamic
principles: reaching inner planets requires reducing angular momentum unless one has a
highly eccentric transfer orbit, which is rare for low thrust, making such transfers less efficient.
Additionally, mass influences GA effectiveness -- Mercury (a mere 5.5\% of Earth's mass) is particularly
unattractive. Although a Venus GA has been demonstrated (e.g. the JUICE mission), its inner position
relative to its departing and arrival body still makes a Venus GA less attractive overall. 

For higher-fidelity methods, which would follow up our RTBA results, selecting optimal sequences in
an early design stage is crucial due to computational costs. A distinct low-$f_s$ group emerges
across all $q$ values, as shown in \cref{figure:results_highfit}, defined by sequences with minimum
$f_s$ within 30\% or 5 km/s of the lowest $f_s$. These sequences maintain consistency despite
variations in seed, $\xi$, and $\chi$. On top of the aforementioned trends in the type of
sub-sequences that are observed in more optimal transfer, this group distinguishes itself by tending
towards numerous Earth flybys to increase orbital energy gradually before transferring to Mars and
subsequently Jupiter. Furthermore, the result aligns well with \cite{fan2021fast}; besides the
omission of EEEMJ in $q = 0.1$, all sequences from \cite{fan2021fast} appear in this optimal group
(in \cref{figure:results_highfit} shown in red), which verifies their optimality, thereby making the
RTBA a reliable way for the first selection of attractive sequences.

\begin{table}[h!]
  \caption{Simulation characteristics per $q$ value. The average run time is rounded to the nearest
  half hour for clarity. The memory used is specifically the RAM usage across all cores, not the
disk space.}
  \label{tab:compchar}
  \begin{center}
    \begin{tabular}[c]{|c|m{2cm}|m{2cm}|}
      \hline
      $q$ [-] & Average run time [hours] &  Average memory utilized [GB] \\
      \hline
      0.1 & \makecell{5} & \makecell{6.80} \\
      0.3 & \makecell{13} & \makecell{7.50} \\
      0.5 & \makecell{21.5} & \makecell{8.15} \\
      \hline
    \end{tabular}
  \end{center}
\end{table}

\cref{tab:compchar} gives an overview of the computational load for various settings. As expected,
the run time for lower $q$ values -- i.e. less sequences evaluated -- is smaller. When considering
the computational complexity, one should look at the efficiency of an algorithm. However, only the
run time is often given in literature, allowing for a limited comparison only. In one example a run
time of 67 hours is needed using 60 Intel XEON E7-4890 @2.8GHz CPUs \cite{englander2017automated}.
By contrast, as can be seen in \cref{tab:compchar}, the 21.5 hours needed by 42 CPUs from DelftBlue
that run on Intel XEON E5-6248R @3.0Ghz CPUs is significantly lower, however, the combinatorial
complexity is different with only three allowed GA candidate planets rather than eight and eight
permissible GAs rather than only three. Another example shows a run time of merely two to three
hours on an Intel Core i5 2.50 GHz CPU for an Earth-Jupiter transfer with three permissible GAs
\cite{fan2022improved}. Furthermore, it can be observed that the computational complexity is
significantly different because the departure date is fixed in \cite{fan2022improved}. However, the
problem definition is incomplete so no true comparison can be made. In general, the run time is very
dependent on the definition of the problem. 

\section{Conclusions} \label{sec:conclusion}

The MGA sequencing problem for low-thrust trajectories can be tackled with an automated
parallel-computing optimization approach by using the RTBA. This strategy consists of a nested-loop
optimization approach that exploits multi-processing capabilities and combines the existing MCTS and
BSS tree-search methods with hodographic shaping. A high-fitness group of sequences was consistently
found for an Earth-Jupiter transfer with up to three GAs, which coincided with those found in
literature and suggests the methodology is robust across sequence lengths. The combination of a
greedy approach with a Monte-Carlo-based search is effective; however, some difficulty arises due to
the limits of using a single PS -- stemming from the greedy nature of the RTBA. By introducing novel
quantities $\xi$ and $\chi$ to determine the 'best' candidate in the RTBA, the convergence was
improved.

%
%

\appendix

\section{Detailed dynamics description} \label{sec:appendix}

As described in \cref{sec:dynamical-description}, two-body motion is used and a continuous
perturbative force is added representing the thrust vector, resulting in
\cref{equation:twobodymotionpert}.

\begin{equation} \label{equation:twobodymotionpert}
  \ddot{\vec{r}} = -G\frac{M}{r^3}\vec{r} + \vec{f}
\end{equation}

Here, $\vec{r}$ is the barycentric satellite position vector, $G$ is the universal gravitational
constant, $M$ is the mass of the attracting body, and $\vec{f}$ is any perturbing acceleration --
defined as $\vec{f}= \frac{\vec{F}}{m_s}$, where $\vec{F}$ is the perturbing force and $m_s$ is the
mass of the spacecraft. 

For the gravity-assist component of the dynamics, a simplified description is used of the
patched-conics approach \cite{conway2010spacecraft}: the gravity-assist maneuver is approximated by
an instantaneous change of momentum rather than a propagated arc within the two-body system of the
GA body. This assumption is acceptable in view of the size of the sphere of influence of the flyby
body, compared to the size of the total transfer. Next, a slightly more detailed explanation is
given of the relevant parameters in a GA because they are used as optimization parameters, described
in \cref{sec:optimapproach}.

%

\begin{figure}
  \begin{center}
    \includegraphics[width=0.45\textwidth]{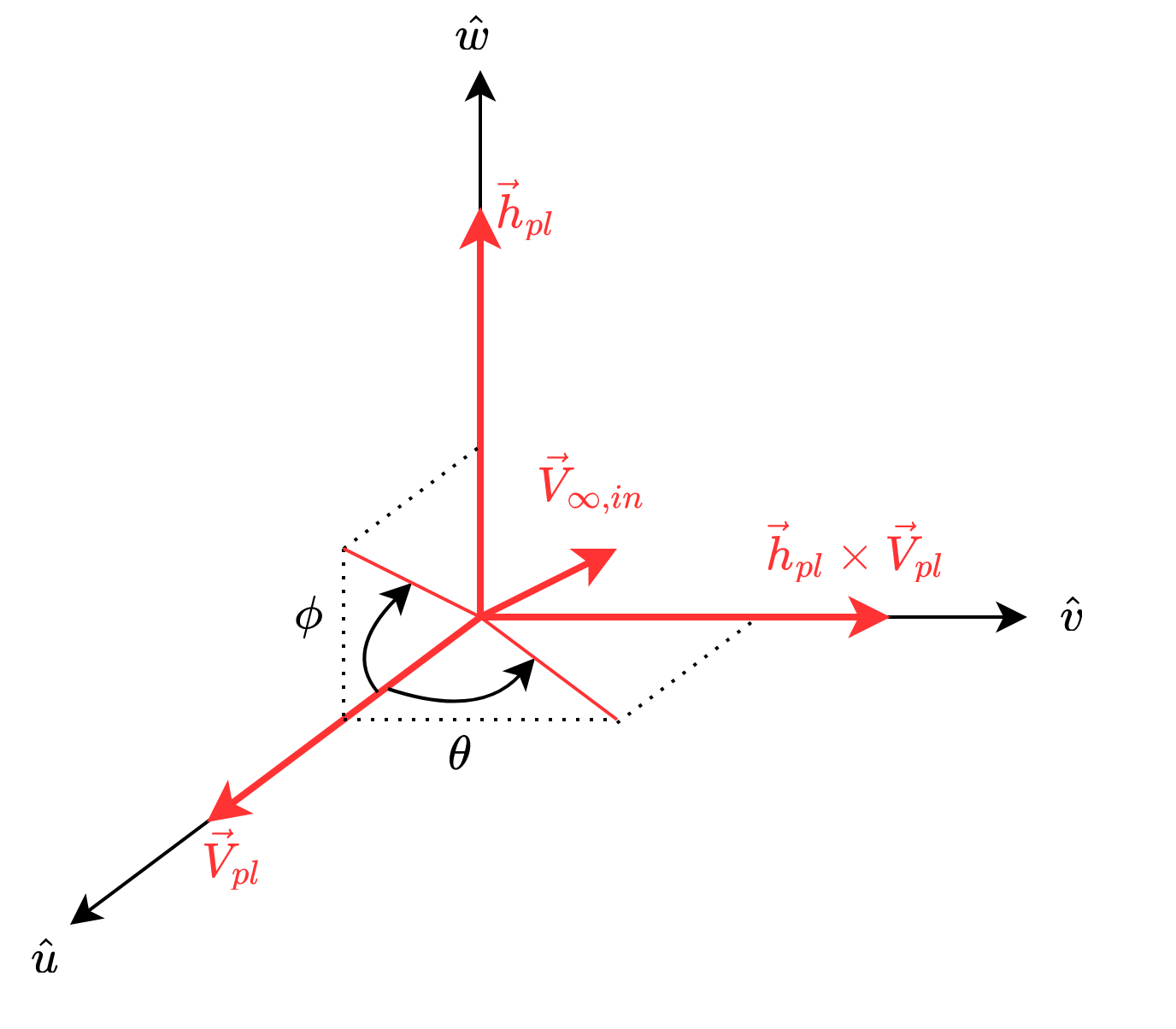}
  \end{center}
  \caption{TNW frame.}
  \label{figure:tnwframe}
\end{figure}

The incoming heliocentric velocity vector is unknown, and is defined using three quantities in
\cref{eq:3.62} \cite{musegaas2013optimization} related to the TNW frame: the norm of the velocity
vector, and $\theta$ and $\phi$, which are the two angular quantities that define the direction of
the incoming hyperbolic velocity vector relative to the planet velocity and angular momentum vector
(see \cref{figure:tnwframe}).

\begin{equation} \label{eq:3.62}
  \vec{V}_{\infty, in} = |\vec{V}_{\infty, in}| (\cos (\phi)\cos (\theta)\hat{u} + \cos (\phi)\sin
  (\theta)\hat{v}+ \sin (\phi)\hat{w})
\end{equation}

Here, $\hat{u}$, $\hat{v}$, and $\hat{w}$ are the orthogonal unit vectors of the TNW frame,
pointing in the planet velocity vector ($\vec{V}_{pl}$), planet angular momentum vector
($\vec{h}_{pl}$), and completing reference frame vector directions ($\vec{h}_{pl} \times
\vec{V}_{pl}$), respectively. To establish an expression for the outgoing hyperbolic velocity
vector, further quantities are needed that link the two gravity-assist related reference frames. 

%
%
%

\begin{figure}
  \begin{center}
    \includegraphics[width=0.45\textwidth]{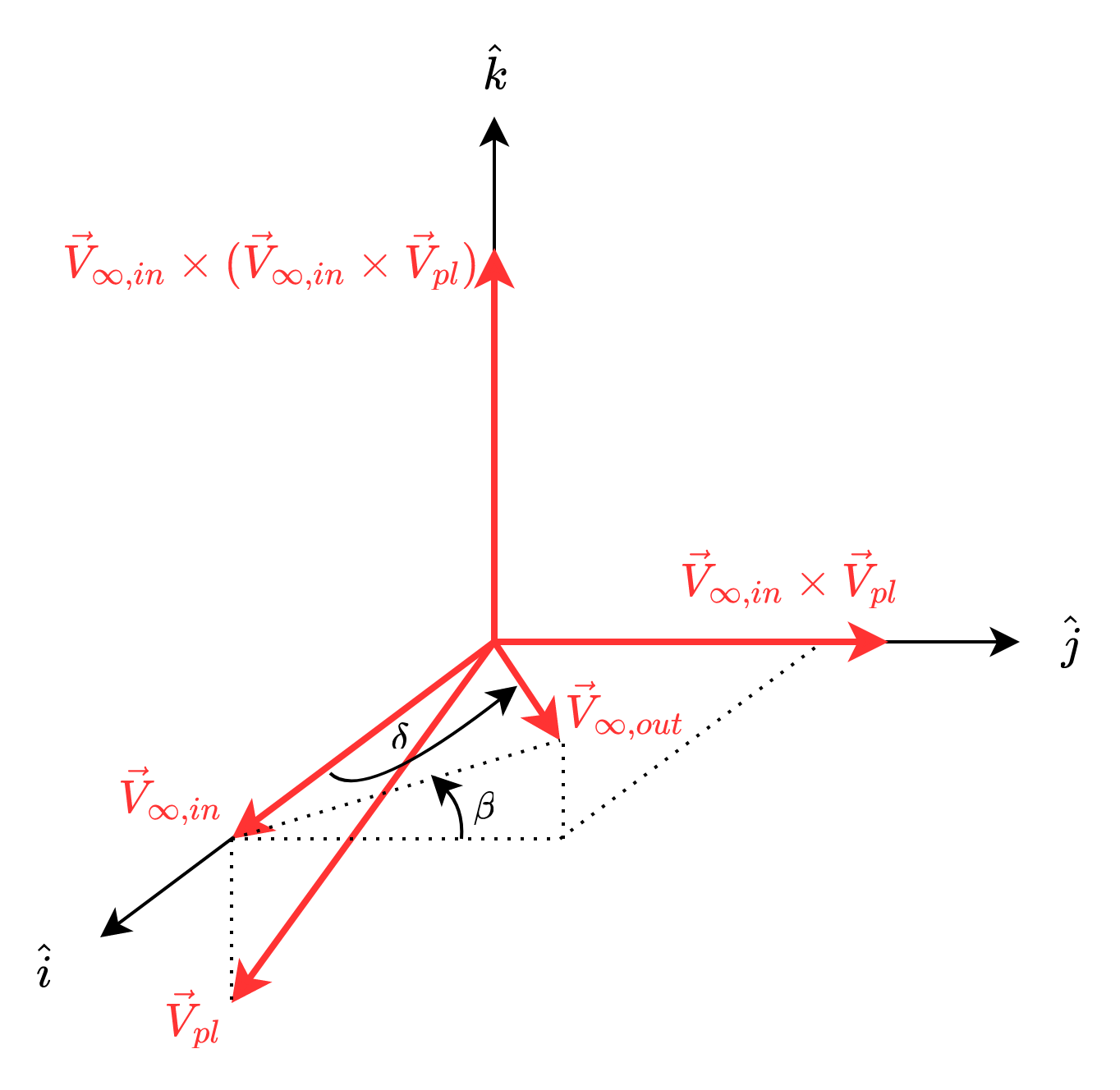}
  \end{center}
  \caption{Local frame.}
  \label{figure:localframe}
\end{figure}

The deflection angle $\delta$ represents the angle between the incoming and outgoing hyperbolic
velocity vector. In addition to the deflection angle, an extra quantity is required to fully define
the nature of the GA maneuver. The outgoing hyperbolic velocity vector is constructed relative to
the incoming hyperbolic velocity vector using not only $\delta$ but also $\beta$, as can be seen in
\cref{figure:localframe,,eq:3.9}, which can be interpreted as the orientation of the plane in which
the GA takes place.

\begin{equation} \label{eq:3.9}
    \hat{V}_{\infty, out} =\cos (\delta) \hat{i}+\cos (\beta) \sin (\delta)
    \hat{j}+\sin (\beta) \sin (\delta) \hat{k}
\end{equation}

In \cref{eq:3.9}, $\hat{V}_{\infty, out}$ is the unit vector for the outgoing hyperbolic velocity
vector, $\beta$ is the plane orientation, $\hat{i} = \frac{\vec{V}_{\infty, in}}{V_{\infty, in}}$,
$\hat{j}=\frac{\hat{i} \wedge \vec{V}^h_{pl}}{|\hat{i} \wedge \vec{V}^h_{pl}|}$, and
$\hat{k}=\hat{i} \wedge \hat{j}$. To get the heliocentric outgoing velocity vector, this outgoing
hyperbolic velocity unit vector is multiplied by the incoming hyperbolic velocity and subsequently
added to the heliocentric planetary velocity vector, as seen in \cref{equation:3.10,,equation:3.11}. 

\begin{equation}
    \vec{V}_{\infty, out} = |V_{\infty, in}|\hat{V}_{\infty, out}
\label{equation:3.10}
\end{equation}

\begin{equation}
    \vec{V}^h_{out} = \vec{V}^h_{pl} + \vec{V}_{\infty, out}
\label{equation:3.11}
\end{equation}





\printbibliography

@article{fan2022improved,
  author = {Fan, Zichen and Huo, Mingying and Quarta, Alessandro A and Mengali,
            Giovanni and Qi, Naiming},
  journal = {Aerospace Science and Technology},
  pages = {107946},
  publisher = {Elsevier},
  title = {Improved Monte Carlo Tree Search-based approach to low-thrust
           multiple gravity-assist trajectory design},
  volume = {130},
  year = {2022},
}

@article{abdelkhalik2012dynamic,
  author = {Abdelkhalik, Ossama and Gad, Ahmed},
  journal = {Journal of Guidance, Control, and Dynamics},
  number = {2},
  pages = {520--529},
  title = {Dynamic-size multiple populations genetic algorithm for
           multigravity-assist trajectory optimization},
  volume = {35},
  year = {2012},
}

@article{englander2017automated,
  author = {Englander, Jacob A and Conway, Bruce A},
  journal = {Journal of Guidance, Control, and Dynamics},
  number = {1},
  pages = {15--27},
  publisher = {American Institute of Aeronautics and Astronautics},
  title = {Automated solution of the low-thrust interplanetary trajectory
           problem},
  volume = {40},
  year = {2017},
}

@article{morante2021survey,
  author = {Morante, David and Sanjurjo Rivo, Manuel and Soler, Manuel},
  journal = {Aerospace},
  number = {3},
  pages = {88},
  publisher = {Multidisciplinary Digital Publishing Institute},
  title = {A survey on low-thrust trajectory optimization approaches},
  volume = {8},
  year = {2021},
}

@inproceedings{hennes2015interplanetary,
  author = {Hennes, Daniel and Izzo, Dario},
  booktitle = {Twenty-Fourth International Joint Conference on Artificial
               Intelligence},
  address = {Buenos Aires, Argentina},
  title = {Interplanetary trajectory planning with Monte Carlo tree search},
  year = {2015},
}

@article{CERIOTTI20101202,
  author = {Matteo Ceriotti and Massimiliano Vasile},
  doi = {https://doi.org/10.1016/j.actaastro.2010.07.001},
  issn = {0094-5765},
  journal = {Acta Astronautica},
  number = {9},
  pages = {1202-1217},
  title = {MGA trajectory planning with an ACO-inspired algorithm},
  url = {https://www.sciencedirect.com/science/article/pii/S0094576510002468},
  volume = {67},
  year = {2010},
}

@phdthesis{ellison2018robust,
  author = {Ellison, Donald Hamilton},
  school = {University of Illinois at Urbana-Champaign},
  title = {Robust preliminary design for multiple gravity assist spacecraft
           trajectories},
  year = {2018},
}

@conference{bellomemulti,
  author = {Bellome, Andrea and Sanchez Cuartielles, Joan Pau and Kemble,
            Stephen and Felicetti, Leonard},
  year = {2021},
  month = {06},
  booktitle = {8th International Conference on Astrodynamics Tools and
               Techniques},
  title = {A multi-fidelity optimization process for complex multiple gravity
           assist trajectory design},
}

@article{gondelach2015hodographic,
  author = {Gondelach, David J and Noomen, Ron},
  journal = {Journal of Spacecraft and Rockets},
  number = {3},
  pages = {728--738},
  publisher = {American Institute of Aeronautics and Astronautics},
  title = {Hodographic-shaping method for low-thrust interplanetary trajectory
           design},
  volume = {52},
  year = {2015},
}

@article{chilan2013automated,
  author = {Chilan, Christian M and Conway, Bruce A},
  journal = {Journal of Guidance, Control, and Dynamics},
  number = {5},
  pages = {1410--1424},
  publisher = {American Institute of Aeronautics and Astronautics},
  title = {Automated design of multiphase space missions using hybrid optimal
           control},
  volume = {36},
  year = {2013},
}

@book{conway2010spacecraft,
  author = {Conway, Bruce A},
  publisher = {Cambridge University Press},
  title = {Spacecraft trajectory optimization},
  volume = {29},
  year = {2010},
}

@article{novak2011improved,
  author = {Novak, Daniel M and Vasile, Massimiliano},
  journal = {Journal of Guidance, Control and Dynamics},
  number = {1},
  pages = {128--147},
  title = {Improved shaping approach to the preliminary design of low-thrust
           trajectories},
  volume = {34},
  year = {2011},
}

@inproceedings{james2017ananalysis,
  author = {James, Steven and Konidaris, George and Rosman, Benjamin},
  booktitle = {Thirty-First AAAI Conference on Artificial Intelligence},
  month = {02},
  number = {1},
  title = {An Analysis of Monte Carlo Tree Search},
  address = {San Francisco, California USA},
  volume = {31},
  year = {2017},
}

@article{wall2009shape,
  author = {Wall, Bradley J and Conway, Bruce A},
  journal = {Journal of Guidance, Control, and Dynamics},
  number = {1},
  pages = {95--101},
  title = {Shape-based approach to low-thrust rendezvous trajectory design},
  volume = {32},
  year = {2009},
}

@article{carnelli2009evolutionary,
  author = {Carnelli, Ian and Dachwald, Bernd and Vasile, Massimiliano},
  journal = {Journal of Guidance, Control, and Dynamics},
  number = {2},
  pages = {616--625},
  title = {Evolutionary neurocontrol: A novel method for low-thrust
           gravity-assist trajectory optimization},
  volume = {32},
  year = {2009},
}

@article{morante2019multi,
  author = {Morante, David and Sanjurjo Rivo, Manuel and Soler, Manuel},
  journal = {Journal of Guidance, Control, and Dynamics},
  number = {3},
  pages = {476--490},
  publisher = {American Institute of Aeronautics and Astronautics},
  title = {Multi-objective low-thrust interplanetary trajectory optimization
           based on generalized logarithmic spirals},
  volume = {42},
  year = {2019},
}

@article{biscani2020parallel,
  author = {Biscani, Francesco and Izzo, Dario},
  journal = {Journal of Open Source Software},
  number = {53},
  pages = {2338},
  title = {A parallel global multiobjective framework for optimization: pagmo},
  volume = {5},
  year = {2020},
}

@article{petropoulos2004shape,
  author = {Petropoulos, Anastassios E and Longuski, James M},
  journal = {Journal of Spacecraft and Rockets},
  number = {5},
  pages = {787--796},
  title = {Shape-based algorithm for the automated design of low-thrust, gravity
           assist trajectories},
  volume = {41},
  year = {2004},
}

@article{vasile2002design,
  author = {Vasile, M and Bernelli-Zazzera, F and Fornasari, N and Masarati, P},
  journal = {Final report of ESA/ESOC study contract},
  number = {00},
  title = {Design of interplanetary and lunar missions combining low thrust and
           gravity assists},
  volume = {14126},
  year = {2002},
}

@article{vasile2015incremental,
  author = {Vasile, Massimiliano and Martin, Juan Manuel Romero and Masi, Luca
            and Minisci, Edmondo and Epenoy, Richard and Martinot, Vincent and
            Baig, Jordi Fontdecaba},
  journal = {Acta Astronautica},
  pages = {407--421},
  publisher = {Elsevier},
  title = {Incremental planning of multi-gravity assist trajectories},
  volume = {115},
  year = {2015},
}

@article{gad2011hidden,
  author = {Gad, Ahmed and Abdelkhalik, Ossama},
  journal = {Journal of Spacecraft and Rockets},
  number = {4},
  pages = {629--641},
  title = {Hidden genes genetic algorithm for multi-gravity-assist trajectories
           optimization},
  volume = {48},
  year = {2011},
}

@inproceedings{sims1999preliminary,
  author = {Sims, Jonathan and Flanagan, Steve},
  title = {Preliminary design of low-thrust interplanetary missions},
  address = {Girdwood, Alaska},
  booktitle = {American Astronautical Society Paper},
  number = {99-338},
  publisher = {AAS/AIAA Astrodynamics Specialist Conference},
  year = {1999},
  
}

@article{strange2002graphical,
  author = {Strange, Nathan J and Longuski, James M},
  journal = {Journal of Spacecraft and Rockets},
  number = {1},
  pages = {9--16},
  title = {Graphical method for gravity-assist trajectory design},
  volume = {39},
  year = {2002},
}

@inproceedings{debban2002design,
  author = {Debban, Theresa and McConaghy, T and Longuski, James},
  booktitle = {AIAA/AAS Astrodynamics Specialist Conference and Exhibit},
  pages = {4729},
  title = {Design and optimization of low-thrust gravity-assist trajectories to
           selected planets},
  address = {Monterey, California},
  year = {2002},
}

@article{crain2000interplanetary,
  author = {Crain, Timothy and Bishop, Robert H and Fowler, Wallace and Rock,
            Kenneth},
  journal = {Journal of Spacecraft and Rockets},
  number = {4},
  pages = {468--474},
  title = {Interplanetary flyby mission optimization using a hybrid global-local
           search method},
  volume = {37},
  year = {2000},
}

@inproceedings{zhang2015analysis,
  author = {Zhang, Jin and Luo, Yazhong and Li, Haiyang and Tang, Guojin},
  booktitle = {2015 IEEE Congress on Evolutionary Computation (CEC)},
  doi = {10.1109/CEC.2015.7256945},
  pages = {596-602},
  title = {Analysis of multiple asteroids rendezvous optimization using genetic
           algorithms},
  address = {Sendai, Japan},
  year = {2015},
}

@article{englander2012automated,
  author = {Englander, Jacob A and Conway, Bruce A and Williams, Trevor},
  journal = {Journal of Guidance, Control, and Dynamics},
  number = {6},
  pages = {1878--1887},
  title = {Automated mission planning via evolutionary algorithms},
  volume = {35},
  year = {2012},
}

@misc{musegaas2013optimization,
  author = {Musegaas, Paul},
  howpublished = {Delft University of Technology \url{https://repository.tudelft.nl/}},
  title = {Optimization of space trajectories including multiple gravity assists
           and deep space maneuvers},
  year = {2013},
  note = {\uppercase{ms}c Thesis},
}

@article{fan2021fast,
  author = {Fan, Zichen and Huo, Mingying and Qi, Ji and Qi, Naiming},
  journal = {Acta Astronautica},
  pages = {233--240},
  publisher = {Elsevier},
  title = {Fast initial design of low-thrust multiple gravity-assist
           three-dimensional trajectories based on the Bezier shape-based method},
  volume = {178},
  year = {2021},
}

@manual{juicemission,
  author={ESA},
  year={n.d.},
  title = {Juice: Jupiter icy Moons Explorer},
  url = {https://www.esa.int/Science\_Exploration/Space\_Science/Juice},
  note = {\uppercase{A}ccessed 2025-05-20},
}

@manual{lisamission,
  author={Tyler, Pat},
  year={n.d.},
  title = {LISA: Laser Interferometer Space Antenna},
  url = {https://lisa.nasa.gov/},
  note = {\uppercase{A}ccessed 2025-05-20},
}

@manual{clippermission,
  author={NASA},
  year={2025},
  title = {Europa Clipper},
  url = {https://science.nasa.gov/mission/europa-clipper/},
  note = {\uppercase{A}ccessed 2025-05-20},
}

@book{pontryagin2018mathematical,
  title={Mathematical theory of optimal processes},
  author={Pontryagin, Lev Semenovich},
  year={2018},
  publisher={Routledge}
}

@inproceedings{rosa2006genetic,
  title={Genetic algorithm and indirect method coupling for low-thrust trajectory optimization},
  author={Rosa Sentinella, Matteo and Casalino, Lorenzo},
  booktitle={42nd AIAA/ASME/SAE/ASEE Joint Propulsion Conference \& Exhibit},
  address = {Sacramento, California},
  pages={4468},
  year={2006}
}

@article{zhang2015low,
  title={Low-thrust minimum-fuel optimization in the circular restricted three-body problem},
  author={Zhang, Chen and Topputo, Francesco and Bernelli-Zazzera, Franco and Zhao, Yu-Shan},
  journal={Journal of Guidance, Control, and Dynamics},
  volume={38},
  number={8},
  pages={1501--1510},
  year={2015},
  publisher={American Institute of Aeronautics and Astronautics}
}

@article{morante2020hybrid,
  title={Hybrid multi-objective orbit-raising optimization with operational constraints},
  author={Morante, David and Sanjurjo-Rivo, Manuel and Soler, Manuel and S{\'a}nchez-P{\'e}rez, Jos{\'e} Manuel},
  journal={Acta Astronautica},
  volume={175},
  pages={447--461},
  year={2020},
  publisher={Elsevier}
}

@article{gao2007near,
  title={Near-optimal very low-thrust earth-orbit transfers and guidance schemes},
  author={Gao, Yang},
  journal={Journal of guidance, control, and dynamics},
  volume={30},
  number={2},
  pages={529--539},
  year={2007}
}

@article{kechichian2009optimal,
  title={Optimal low-thrust transfer in general circular orbit using analytic averaging of the system dynamics},
  author={K{\'e}chichian, Jean A},
  journal={The Journal of the Astronautical Sciences},
  volume={57},
  number={1},
  pages={369--392},
  year={2009},
  publisher={Springer}
}

@manual{eclipj2000,
  author={NASA},
  year={n.d.},
  title = {Reference Frames},
  url = {https://naif.jpl.nasa.gov/pub/naif/toolkit\_docs/C/req/frames.html},
  note = {\uppercase{A}ccessed 2025-06-23},
}







\end{document}